\documentclass[10pt,letterpaper]{article}
\usepackage[top=0.85in,left=2.75in,footskip=0.75in]{geometry}

\usepackage{amsmath,amssymb}

\usepackage{changepage}

\usepackage{textcomp,marvosym}

\usepackage{cite}

\usepackage{nameref,hyperref}
\usepackage{tabularx}

\usepackage[right]{lineno}

\usepackage[nopatch=eqnum]{microtype}
\DisableLigatures[f]{encoding = *, family = * }

\usepackage[table]{xcolor}

\usepackage{array}

\newcolumntype{+}{!{\vrule width 2pt}}

\newlength\savedwidth

\raggedright
\usepackage[aboveskip=1pt,labelfont=bf,labelsep=period, justification=raggedright ,singlelinecheck=off]{caption}

\makeatletter
\renewcommand{\@biblabel}[1]{\quad#1.}
\makeatother

\usepackage{lastpage,fancyhdr,graphicx}
\usepackage{epstopdf}
\usepackage{todonotes}
\reversemarginpar
\fancyheadoffset[L]{2.25in}
\fancyfootoffset[L]{2.25in}
\newcommand\fay[1]{\textcolor{violet}{#1}}

\begin{document}
\vspace*{0.2in}

\begin{flushleft}
{\Large
\textbf\newline{Boolean Algebra–Driven Sepsis Diagnosis } 
}
\newline
\\
Marcus Weber\textsuperscript{1*},
Kai Kappert\textsuperscript{2},
Marco Reidelbach\textsuperscript{1}, 
Ambros Gleixner\textsuperscript{3},
Konstantin Fackeldey\textsuperscript{1,4*},
Wolfgang Bauer\textsuperscript{5}
\\
\bigskip
\textbf{1} Zuse Institute Berlin, Takustraße 7, D-14195 Berlin, Germany\\
\textbf{2}  Charité - Universitätsmedizin Berlin, corporate member of Freie Universität Berlin and Humboldt-Universität zu Berlin, Institute of Diagnostic Laboratory Medicine, Clinical Chemistry and Pathobiochemistry, Augustenburger Platz 1, D-13353 Berlin, Germany. \\
\textbf{3} HTW Berlin – Hochschule für Technik und Wirtschaft Berlin, D-10313 Berlin, Germany  \\
\textbf{4} Technische Universität Berlin, Institut für Mathematik, Straße des 17. Juni 136, D-10623 Berlin, Germany  \\
\textbf{5} Charité - Universitätsmedizin Berlin, corporate member of Freie Universität Berlin and Humboldt-Universität zu Berlin, Department of Emergency Medicine, Hindenburgdamm 30, D-12203 Berlin, Germany

\bigskip

%
%





* weber@zib.de, fackeldey@math.tu-berlin.de (corresponding authors)

\end{flushleft}

\section*{Abstract}
Sepsis remains a diagnostic challenge due to its heterogeneous molecular signatures and complex immune responses. In this study, we develop a logical data analysis framework based on Boolean polynomial rings. This method constructs an ideal $\mathcal{I}$ of selection criteria that isolate empty subsets of previously analyzed patient data. This approach enables the derivation of interpretable classification rules based on biomarker profiles. We demonstrate that logical data analysis identifies distinct logical patterns for positive and negative sepsis classification. For instance, elevated levels of GLP-1 and MyD88 are associated with septic states in our dataset, whereas high TRAIL and low MyD88 concentrations may suggest a non-septic condition. Importantly, a new way to integrate expert knowledge to filter out potential overfitting or dataset-specific artifacts is shown. Our findings highlight the utility of logics in generating transparent, biologically plausible rules for a data-based and expert-based understanding of sepsis. Moreover, we show how data analysis can benefit from algebraic structures. 


\section*{Introduction}

Sepsis is defined as life-threatening organ dysfunction caused by a dysregulated immune response of the host to an infection \cite{Singer2016}. It is a life-threatening disease with a high incidence. Rudd, Johnson, Agesa et al. used multicausal cause of death statistics for 2017 to determine a global incidence of 677.5 per 100,000 person-years \cite{Rudd2020}. Based on these data,
approximately 48.9 million patients worldwide developed sepsis in 2017, of whom 11 million died. Approximately one in five people worldwide died in connection with sepsis.
 
The pathophysiology of sepsis is not yet fully understood, and the clinical manifestation is inconsistent and often nonspecific in the initial phase. However, early identification and rapid initiation of therapy are crucial for the prognosis of patients with sepsis, as a delay in adequate therapy is directly related to mortality \cite{Seymour2017}. The diagnosis of sepsis is difficult because it arises from a complex interaction between the infection and the host response, which varies widely between individuals. Factors such as the type and severity of the infection, as well as host characteristics including age, comorbidities, and immune status, determine whether an infection progresses to sepsis \cite{Tuerxun_2023}. Compared to other life-threatening diseases, such as myocardial infarction, for which precise and multimodal diagnostics (biomarkers, ECG, imaging) are available, there is currently no diagnostic test that can reliably identify and predict patients with sepsis or potentially septic disease progression.

This article introduces an innovative approach for developing  logical statements based on initial numerical and clinical measures, more concretely: assessing concentrations of proteins in blood serum. In this article we try to answer the question whether it is possible to derive logical statements from blood serum concentration measurements of proteins to support the diagnosis of sepsis. Early attempts to solve this problem can be found in 2014 \cite{Gultepe} where the authors use classification methods for their decision support system. The authors of \cite{Pepic2021} summarize existing methods for early detection of sepsis using artificial intelligence. Modern medical diagnostics are increasingly shaped by mathematical models and technical systems. AI-supported methods and probabilistic modeling are used, particularly for complex clinical conditions such as sepsis. Arriaga-Pizano, Gonzalez-Olvera, Ferat-Osorio, et al. \cite{arriaga2021sepsisNN} demonstrated that neural networks can enable sepsis detection based on routine clinical variables. Yao, Yang, Xu, et al. \cite{yao2021probabilistic} developed a probabilistic rule model with a focus on interpretability, while Zhang, Zhong, Cheng, et al. \cite{zhang2025sepsisPlatform} presented a real-time platform with dynamic time series features and a visualization. However, all of these models focus on patients who are already in the intensive care unit.

At the same time, Boolean algebra is establishing itself as a structuring tool for formalizing medical decision-making processes. By converting complex decision processes into logical rules, it can also serve as a traceable support.  In biochemistry, logic-based methods model cellular signaling networks \cite{logicSignaling2010}, while Macauley and Youngs \cite{macauley2020algebraic} emphasize algebraic approaches as key concepts in biological systems analysis. Zugmaier et al. \cite{zugmeier2023myeloid} demonstrated the application of Boolean logic to the classification of hematological neoplasms, and other studies use logical operators to improve laboratory-based diagnostics \cite{boolean2024lab} and for computer-aided differential diagnosis \cite{boolean2025diagnosis}. The combination of data-driven AI and mathematical-logical modeling (outside the field of sepsis) provides a robust foundation for transparent, automated, and reproducible diagnostic systems. This is supported by reviews \cite{review2019sepsisAI} and clinical evaluations \cite{gupta2024computational}, which demonstrate the bridge between technology and patient care. In particular when it comes to establishing decision trees \cite{costa_recent_2023} as a diagnostic tool. The structured, rule-based nature of these trees allows the development of suitable and transparent diagnostic procedures for everyday clinical practice.  

The analyses referenced above derive the evidential basis for their claims directly from the content of the datasets themselves. Rather than relying on external contextual knowledge or on established experts' assumptions, these approaches emphasize internal consistency, pattern recognition, and inference grounded in the data. In most cases, the conclusions are supported through statistical reasoning -- employing techniques such as correlation analysis, probabilistic modeling, or algorithmic classification -- to ensure that the interpretations are both reproducible and  quantitatively (i.e., statistically) robust.

However, certain principles that play a central role in the analysis of datasets within the humanities -- especially using contextual knowledge not contained in the datasets themselves for the interpretation of data -- are also of great interest in the diagnosis of diseases. This shared analytical foundation forms the bridge between the two domains.  In this work, the authors aim to establish a conceptual and methodological bridge between two distinct yet innovative studies: the archaeological investigation of structures in ancient settlements by using logical data analysis \cite{Klasse2025}, and the biomedical application of composite algorithmic analysis for sepsis diagnostics \cite{BAUER2025106599}. While the former explores logical data analysis within the context of the Workmen’s Village in Tell el-Amarna, Egypt, the latter study applies statistical principles to evaluate biomarkers for sepsis in clinical settings. This study will provide the dataset for our investigations.

\section*{Materials}

The design of the studies used for protein abundance in patients suspected for infection and sepsis, have previously been described elsewhere \cite{Bauer2021b, Diehl-Wiesenecker2024, 10.1093/ofid/ofac437, BAUER2025106599}. Briefly, within study enrollment, patients for whom an acute infection was suspected (as leading differential diagnoses) were identified by trained emergency department (ED) staff. Between February 2019 and August 2022 patients were recruited at the point of presentation to the ED, and serum samples were collected simultaneously with the earliest routine sampling, and proteins were measured as described in \cite{BAUER2025106599}.

\subsection*{Parameters of the dataset}

In the mentioned study, nine different proteins were focused on as potential biomarkers: ENA-78, Fractalkine, GLP-1, Leptin, MMP-8, MyD88, PD-L1, Pentraxin-3 and TRAIL. The pre-selection of potential biomarkers was based on expert knowledge of mechanisms related to infection and inflammation, and on hypotheses about sepsis pathophysiology. The selected proteins were considered promising because previous smaller studies suggested that they might be able to detect the transition from infection to sepsis, but they have not yet been sufficiently investigated in a bigger clinical investigation.
Using a retrospective definition of sepsis (not based on the mentioned biomarkers but on infection states and on clinical outcomes), the ($n=390$) measurements of the study (rows of the dataset) have finally been classified as ``sepsis'' (137) or ``non-sepsis'' (253) measurements. The goal of logical data analysis is to find true logical statements about the concentration of the nine proteins in blood serum and the corresponding sepsis classification. To perform logical data analysis like in \cite{Klasse2025} we need to define Boolean variables. 

\subsection*{Turning numerical data into Boolean data}

Imagine a high concentration of one of the nine proteins would sufficiently discriminate between sepsis and non-sepsis patients. In this case, one could sort the patients according to the concentration of this protein and the $137$ highest concentrations would indicate sepsis, whereas the $253$ lowest concentrations would indicate non-sepsis.  This is the heuristic starting point to classify concentrations of proteins as ``high'' (1) and ``not high'' (0): the thresholds for the protein concentrations are selected in such a way that always in exactly $137$ cases the classification is ``high concentration''. The resulting thresholds are shown in Tab.~\ref{tab:extremes}. Furthermore, this table shows which Boolean variables ($E,F,G,L,M,y,P,x,T$) are used in the algebraic algorithms and in this article. These variables have two possible values (0 or 1). Additionally the classification variable for sepsis is given by $s\in\{0,1\}$. The table of $\{0,1\}$-patterns of the patients was the starting point of the presented analysis.  

\begin{table}[th]
\centering
\begin{tabular}{|l||c|c|}
\hline
protein& threshold\cr
\hline
\hline
(E)NA-78 & $>$ 786.51\cr
\hline
(F)ractalkine & $>$ 	17307.21\cr
\hline
(G)LP-1 & $>$ 22.71\cr
\hline
(L)eptin & $>$ 20818.6\cr
\hline
(M)MP-8 & $>$ 59704.08\cr
\hline
 M(y)D88 & $>$ 51.88\cr
\hline
(P)D-L1 & $>$ 844.84	\cr
\hline
 Pentra(x)in-3 & $>$ 28200.36\cr
\hline
(T)RAIL & $>$ 92.8 \cr
\hline
\end{tabular}\caption{\label{tab:extremes}First column: The 9 selected proteins. Second column: The resulting threshold values for classification as ``high concentration''.}

\end{table}

\section*{Theory}

The aim of logical data analysis is to identify true logical statements about the patients and their sepsis classification. In principle,  statements like ``If the patient has a high concentration of protein A and a low concentration of protein B, then he/she does not have sepsis'' are looked for.  

\subsection*{True logical statements and selection criteria}

Mathematically speaking, a true logical statement can be understood as an equation of selection criteria \cite{Klasse2025}. The expression $A\cdot (B+1)$ with the Boolean variables $A$ and $B$ is true if and only if $A=1$ and $B=0$. In this expression multiplication ($\cdot$) represents the logical operation ``... and ...'', whereas summation ($+$) represents the logical operation ``either... or ...''\fay{\footnote{Also known as XOR operation $\oplus$ in the Boolean literature.}}. The equation of selection criteria
\[
A\cdot(B+1) = A\cdot(B+1)\cdot (s+1)
\]
means that selecting all patients with a high concentration of $A$ and a low concentration of $B$ is the same like selecting all patients with a high concentration of $A$ and a low concentration of $B$ and not having sepsis. This equation corresponds to the true logical statement ``If the patient has a high concentration of protein A and a low concentration of protein B, then he/she does not have sepsis''. One can do Boolean ring transformations with equations like these. By factoring out we get:
\[
AB+A = ABs + As + AB + A.
\]
Adding $AB+A$ on both sides and using $(AB+A)+(AB+A)=0$ (``either $(AB+A)$ or $(AB+A)$'' is always false)  leads to:
\[
0 = A(B+1)s.
\]
The expression $0=A(B+1)s$ means that the statement ``If the patient has a high concentration of protein A and a low concentration of protein B, then he/she does not have sepsis'' is equivalent to say that there is no patient with a high concentration of $A$, a low concentration of $B$, and having sepsis $s$. 

Let $S_1$ and $S_2$ be two Boolean expressions representing two selection criteria. Then the equality $S_1=S_2$ means that both selection criteria have identical truth values for each patient. This expression can be equivalently transformed into $S_1+S_2=0$, which states that there do not exist patients who differ in these truth values.

\vspace*{0.3cm}
\fbox{\begin{minipage}{0.94\textwidth}\em A true logical statement which is expressed by an equation of Boolean expressions $(S_1=S_2)$ with identical truth values for each observation is equivalent to a Boolean expression $(S_1+S_2=0)$ which is always false. $S_1+S_2$ applied as a selection criterion does not select an observed measurement. \end{minipage}} \vspace*{10pt}  


\subsection*{Number of empty selection criteria}

Searching for true logical statements is to identify those selection criteria which do not ``find'' patients. 
In our example we have found 194 different $\{0,1\}$ patterns among the 390 measurements which correspond to observed combinations of values of the $10$ Boolean variables. How many selection criteria do not ``find'' observed patterns? 

Having 10 binary digits at hand, i.e., 9 proteins plus ($s$)epsis, we can in principle create $2^{10}=1024$ different patterns. Out of these possible patterns, we only have observed $194$. We have not observed $830$ possible patterns. How many selection criteria exist which select a {\bf subset} of these $830$ non-observed patterns? It is $2^{830}$, because every non-observed pattern can either be part of such a subset or not. This means, our given data provides the possibility to formulate $2^{830}$ ``empty'' selection criteria leading to true logical statements about the relationship of the $10$ Boolean variables. This number is much bigger than the number of atoms in the universe. 

It is impossible for an expert to handle all true logical statements about  the 10 parameters of the investigated patients. We will have to find a way to pick a relevant subset of these true logical statements. 

\subsection*{Algebraic structure of all empty selection criteria}

We search for selection criteria, which do not ``find'' patients. The set of all possible selection criteria that chose only patterns (combinations of measurements) that have not been observed is denoted as $\mathcal{I}$.  The two polynomials 
\begin{equation}\label{eq:selection}
FTs(y+1), \quad (F+1)xTs,
\end{equation}
are elements of the set $\mathcal{I}$, which means that applying these selection criteria to the given dataset with thresholds defined in Tab.~\ref{tab:extremes} leads to $FTs(y+1)=0$ and $(F+1)xTs=0$, respectively, for all patients.
The criteria in Eq.~\ref{eq:selection} are two different selection criteria, which select a subset of possible patterns, but these patterns have not been observed. There has been no patient with an increased concentration of Fractalkine and TRAIL, a low concentration of MyD88 who has sepsis. There has been no patient with low concentration of Fractalkine and high concentrations of Pentraxin-3 and TRAIL having sepsis. The true logical statements we can derive from these non-observations Eq. \ref{eq:selection} are: 
\begin{itemize}
    \item[R1)] If a patient has a high concentration of TRAIL and Fractalkine, but a low concentration of MyD88, then he/she does not have sepsis.
    \item[R2)] If a patient has a high concentration of TRAIL and Pentraxin-3, but a low concentration of Fractalkine, then he/she does not have sepsis. 
\end{itemize}

If we have found a selection criterion that does not ``find'' a patient, then restricting this statement also does not ``find'' a patient and thus leads to a further true logical statement. Mathematically speaking, if we multiply a suitable selection criterion with {\em any} Boolean polynomial including our 10 variables, it leads to another selection criterion in the set $\mathcal{I}$. For example,
\begin{equation}\label{eq:selection2}
FTs(y+1) \cdot E 
\end{equation}
is a restriction of $FTs(y+1)$ which leads to the
\begin{itemize}
\item[R3)] If a patient has a high concentration of TRAIL and Fractalkine {\em and} ENA-78, but a low concentration of MyD88, then he/she does not have sepsis. 
\end{itemize}

If we have found two suitable selection criteria $S_1$ and $S_2$ which select an empty set of patients, then the selection criterion $S_1 + S_2$ (meaning either $S_1$ or $S_2$) also selects the empty set. In our example we can add the two expressions from Eq. \ref{eq:selection} and arrive at
\[
FTs(y+1) + (F+1)xTs = FTsy+FTs+FxTs+xTs=0.
\]

If $\mathcal{I}$ is the set of selection criteria which do not ``find'' patients, then the discussed properties mean that for any statement $S_1\in{\mathcal{I}}$ and any Boolean polynomial $\pi$, we have $S_1\cdot \pi\in{\mathcal{I}}$. We also discussed that for two expressions $S_1,S_2\in{\mathcal{I}}$ we have $S_1+S_2\in {\mathcal{I}}$. These two findings simply mean that $\mathcal{I}$ is an {\em ideal} inside the set of Boolean polynomials. In our case, the ideal consists of $2^{830}$ elements (including the $0$ expression). The strong algebraic structure makes it easier to handle this big number. The idea to use Boolean rings for the analysis of non-mathematical data sources has been developed in a project together with humanities \cite{weber2020complexity, weber2022mathematics}.

\subsection*{Gröbner Basis of the ideal}
It is an exercise in Commutative Algebra textbooks (e.g., \cite{At1994}) to show that every finitely generated ideal in a Boolean ring is a {\em principal ideal}. The ideal $\mathcal{I}$ is a principal ideal which means that it is generated by just {\em one} selection criterion. In our context, this mathematical finding is also explainable in a different way, because if we have found the selection criterion $\sigma\in {\mathcal{I}}$ which exactly selects all patterns which have not been observed, then the ideal $\mathcal{I}$ just consists of all restrictions $\pi \cdot \sigma$ of this selection criterion, where $\pi$ is a Boolean polynomial of our $10$ variables. $\sigma$ generates $\mathcal{I}$. 

The advantage of knowing the algebraic structure of $\mathcal{I}$: We can easily construct the selection criterion $\sigma$. If $S_1$ selects {\em one existing/observed} pattern and $S_2$ selects {\em one  different} observed pattern, then $S_1+S_2$ selects exactly the union set of these two patterns. Thus, the sum of all different Boolean polynomials constructed according to the last row of Table \ref{tab:pattern} provides the selection criterion $\overline{\sigma}$ which exactly selects all existing/observed patterns. Therefore, the complement $\sigma=1+\overline{\sigma}$ is our searched for generator of the ideal $\mathcal{I}$. Thus, the ideal $\mathcal{I}$ is constructable in computer algebra systems \cite{weber2022coding}. The disadvantage of $\sigma$: If we would write down the generator $\sigma$ of the sepsis data set, it would be a polynomial which fills many pages and is still not manageable.  
\begin{table}[h]
\centering
\small
\begin{tabular}{|p{10pt}p{15pt}p{10pt}p{10pt}p{14pt}p{10pt}p{10pt}p{14pt}p{7pt}p{10pt}|}
\hline
E& F& G& L& M & y& P& x& T& s\cr
\hline
\hline
 1877	&90658	&22	&561	& 60876&	81&	479	&52470	& 408 &yes
\cr
\hline
1&	1&	0&	0&	1&	1&	0&	1&	1&  1\cr
\hline
E & F & G+1 & L+1 & M & y & P+1 & x & T&  s\cr
\hline
\end{tabular}
\caption{\label{tab:pattern}How clinical data is turned into a selection criterion. First row: The ``names'' of the parameters (the variables). Second row: The rounded measurements for patient No. 135. Third row: According to Table \ref{tab:extremes} we transform these measurements into categories $0$ and $1$. This row defines the ``pattern'' of patient No. 135. Fourth row: For building a selection criterion we have to check whether the variable has value $0$ or $1$. The resulting selection criterion is: $E F (G+1) (L+1) M y (P+1) x T s$ . }
\end{table}

Thus, we search for a different basis of this ideal. Looking for a Gr\"obner Basis of the ideal $\mathcal{I}$ has been proposed in \cite{Klasse2025}. Such a basis can be computed using the Buchberger Algorithm applied to the ideal generated by $\sigma$ \cite{Bu70, Be1993}.  There is a good reason for this choice: Looking at Eq. \ref{eq:selection} and Eq. \ref{eq:selection2}, we see that logical restrictions of ``simple polynomials'' are expressed by multiplications. The leading monomial of the polynomial $FTsy+FTs$ in Eq. \ref{eq:selection} is a divider of the leading monomial of the polynomial $EFTsy+EFTs$ in Eq. \ref{eq:selection2}. If a selection criterion $S$ is element of the ideal $\mathcal{I}$, then there should be one element in the chosen basis of $\mathcal{I}$ having a leading monomial which is a divider of the leading monomial of $S$ (to be less ``restrictive'').  The Gr\"obner Basis meets these ``divider requirements''. It is also denoted as standard basis in algebra \cite{Bu70, Hi64}.

\section*{Methods}

Now it has to be described how to derive true logical statements about sepsis out of the elements of the Gröbner Basis of $\mathcal{I}$. The selection criteria which include the variable $s$ provide such logical rules. The general algebraic form of such selection criteria is $S_1\cdot s + S_2$, where $S_1$ and $S_2$ are suitable selection criteria (not necessarily elements of $\mathcal{I}$). 

{\bf Positive Rules.} Imagine $S_2\not\in{\mathcal{I}}$, i.e., $S_2$ selects a subset of patients. The requirement is $S_1\cdot s+S_2\in{\mathcal{I}}$. We can conclude that the selection criterion $S_1\cdot s$ must select the same subset of patients like $S_2$, because it has to ``cancel out with $S_2$'', i.e. $S_1\cdot s=S_2\not\in \mathcal{I}$ if applied to the dataset. This means that the patients in this special subset must have sepsis. 
To give an example: $GMys + Gyxs + GMy + Gyx$ is an element of the Gröbner Basis of $\mathcal{I}$ with $S_1=GMy + Gyx$ and $S_2=GMy + Gyx$.  $S_2$ selects 22 patients and therefore $S_2\not\in \mathcal{I}$. This means, that the 22 patients selected by $S_2=Gy(M+x)$ have sepsis. Patients who have a high concentration of GLP-1 and of MyD88 and either a high concentration of MMP-8 or Pentraxin-3 have sepsis.

{\bf Negative Rules.} Imagine $S_1\cdot(1+S_2)\not\in\mathcal{I}$, i.e., the intersection of the set of patients selected by $S_1\not\in {\mathcal{I}}$ and the complement of the set selected by $S_2$ with $(1+S_2)\not\in{\mathcal{I}}$ is not empty.  In this case, $S_1\cdot s+S_2$ can only select an empty set of patients, if the patients selected by $S_1\cdot(1+S_2)$ do not have sepsis.  To give an example: $FyTs + FTs$ is an element of the Gröbner Basis of $\mathcal{I}$ with $S_1=FyT+FT$ and $S_2=0$, i.e., with $S_1(S_2+1)=FT(y+1)$. This selection criterion selects 24 patients who do not have sepsis. Patients with a high concentration of Fractalkine and TRAIL and a low concentration of MyD88 do not have sepsis. 

\subsection*{Context knowledge and relevance}
In our case study, the Gröbner Basis of $\mathcal{I}$ has more than 200 elements leading to many positive and negative true logical statements about the patients.  Not all of them are regarded as relevant for mainly two different reasons: (i) they only apply to a very small set of patients and can not be justified without further experiments or (ii) they are overfitted with regard to the data. As an example for (i): In a later step of our algorithmic approach we will find a negative rule saying that all patients with high concentration of Pentraxin-3 and TRAIL having a low concentration of Fractalkine  do not have sepsis. However, this selection criterion only selects 7 patients. It may be a ``random'' observation and there might be reasonable doubts that a high concentration of Pentraxin-3 is really a good indicator for excluding sepsis. As an example for (ii): Also in a later step of the algorithm, we find a negative rule applicable to 75 patients and based on the selection criterion $T(EL + EM + EP + Ex + E + FP + GL + GM + Gy +$ $ LM + L + MP + Mx + M + y + Px)$. However, this selection criterion is more of a description of the properties of these 75 patients rather than really an ``insight''.  In principle, we look for ``insights'' among the true logical statements, which is a category that is not only related to the data itself.  

\vspace*{0.3cm}
\fbox{\begin{minipage}{0.94\textwidth}\em Not every true logical statement about sepsis is an insight. The relevance of true logical statements can only be accessed by using experts' or external knowledge not necessarily  provided by the given dataset itself.\end{minipage}} \vspace*{10pt}  

In our algorithmic approach, we will select a subset of polynomials of the Gröbner Basis of $\mathcal{I}$ which are denoted as ``relevant'' and will generate another ideal $\mathcal{J}$ using these polynomials. 

\subsection*{Elimination ideal}
Selecting relevant polynomials from the Gröbner Basis as generator of a subideal ${\mathcal{J}}\subset{\mathcal{I}}$ has another advantage in searching for further relevant true logical statements about the patients.  Polynomials in a Gröbner Basis of $\mathcal{I}$ often start with ``simple polynomials'' and then become increasingly complex (sums of many monomials with higher degree), because they must encode the full structure of the ideal under a chosen monomial ordering. Elimination ideals offer a powerful way to simplify this complexity: Taking simple relevant polynomials as generator of a further ideal $\mathcal{J}$, we can simplify the more complex polynomials of the Gröbner Basis of $\mathcal{I}$ by computing the remainder of the Gröbner polynomials with regard to $\mathcal{J}$. The ``simpler'' remainder polynomials are still elements of $\mathcal{I}$ and can also be used to derive relevant positive and negative true logical statements from them. This procedure of adding polynomials from $\mathcal{I}$ into $\mathcal{J}$ will, therefore, become an 
(inner) iterative loop of our algorithm. 

\subsection*{Exceptions}

The selection criterion $T(y + 1)(L + M)$ selects 45 patients in the given dataset. Except for patient No. 2237, all of the selected patients do not have sepsis. Thus, the negative rule ``patients with $T(y+1)(L+M)$ do not have sepsis'' would be valid, if we exclude patient No. 2237 from our list of observations. By saying that patient No. 2237 is an exception, we make his/her observation a non-observation. This creates new elements in $\mathcal{I}$, i.e., new non-ob\-ser\-va\-tions and also new true logical statements. Cicero's phrase ``exceptio probat regulam'' -- often translated as ``the exception proves the rule'' -- is frequently misunderstood. What Cicero meant was that the existence of an exception implies the presence of a general rule. For example, if a sign says ``No parking on Sundays,'' it suggests that parking is allowed on other days. The exception (Sunday) confirms that a rule (parking is generally permitted) exists. In Cicero’s legal context, this logic was used to infer general principles from specific exclusions. 

The procedure of adding polynomials which select ``exceptional patients'' to the generator of $\mathcal{I}$ to eliminate these observations will, therefore, become an outer iterative loop of our algorithm. The dataset has $n=390$ observations with 194 different patterns. This means that many patterns are only observed once or twice. Not all of these rare observations are exceptions. 

\vspace*{0.3cm}
\fbox{\begin{minipage}{0.94\textwidth}\em Not every rare observation is an exception. The relevance of observations can only be accessed by using experts' insights with regard to the study.\end{minipage}} \vspace*{10pt}  

Also in statistical procedures in which logical rules are claimed on the basis of ``significantly'' frequent observations, exceptions are defined by listing those patients to whom the chosen rules do not apply. ``The rules prove the exceptions.'' 

\subsection*{Algorithm}

Here we will quickly describe the steps of the algorithmic approach to find relevant true logical statements about sepsis based on the given dataset and based on the heuristic choice of thresholds given in Tab.~\ref{tab:extremes}.

\subsubsection*{Step 1: generating the Gröbner Basis of ideal $\mathcal{I}$}

If the generator $\sigma$ of $\mathcal{I}$ is computed, we can generate a Gröbner Basis of this ideal. In the course of the algorithm (in Step 4), polyomials will be determined which will also become further generators and extend the ideal $\mathcal{I} $.

\subsubsection*{Step 2: generating the Gröbner Basis of ideal $\mathcal{J}$}

In the course of the algorithm (in Step 3), polynomials of the Gröbner Basis of $\mathcal{I}$ will be selected and denoted as ``relevant''. These polynomials generate a further ideal $\mathcal{J}$ with its Gröbner Basis. We start with an empty ideal ${\mathcal{J}}=\{0\}$. 

\subsubsection*{Step 3: inner loop to identify ``insights''}

Transform the remainder polynomials of the Gröbner Basis of $\mathcal{I}$ with regard to the ideal $\mathcal{J}$ into positive and negative rules about sepsis. Then go through this list of rules and figure out, which ones are relevant. The respective polynomials of the Gröbner Basis will be added to the list of generators which generate the ideal $\mathcal{J}$. Continue with Step 2, if relevant rules have been found. Otherwise go to Step 4.

\subsubsection*{Step 4: outer loop to identify ``exceptions''}

Go through every observed pattern of the dataset and determine the selection criteria constructed according to Tab.~\ref{tab:pattern}. For every such polynomial compute the remainder with regard to ideal $\mathcal{I}$. Using these remainders one decides whether a rarely observed pattern is really an exception. All remainders, which belong to ``exceptions'' are collected and added to the list of generators of ideal $\mathcal{I}$. If exceptions have been identified, then continue with Step 1. Otherwise the algorithm stops here.  

\section*{Results}

We applied the above algorithm to the given dataset. After three cycles of the outer loop, the algorithm terminated. The polynomials which have been added to list of generators of the ideals $\mathcal{I}$ and $\mathcal{J}$, respectively, in each iteration cycle are shown in Tab.~\ref{tab:results}. In total, we identified 3 patients which we denoted as ``exceptions''. 
\begin{table}[ht]
\centering
\begin{tabular}{|l|p{4cm}|p{4cm}|l|}
\hline
cycle No. & polynomials to $\mathcal{J}$ (inner) & polynomials to $\mathcal{I}$ (outer) & exceptions\cr
\hline
1 & $FyTs + FTs,$ $ GMys + Gyxs + GMy + Gyx$ & $LyPxTs + LPxTs$ & No. 2237\cr
\hline
2 & $LyTs + MyTs + LTs + MTs,$ $ELys + EMys + ELs + EMs$ & $yxTs + xTs$ & No. 127, 545\cr
\hline
3 & $EyTs + MyTs+ yPTs + ETs + MTs + PTs$ & & \cr
\hline
\end{tabular}
\caption{\label{tab:results} In each outer iteration loop (in our example 3 cycles) of the algorithm, polynomials from the Gröbner Basis of $\mathcal{I}$ are added to the list of generating polynomials of the ideal $\mathcal{J}$. These polynomials are listed in the second column. After finalizing the inner loop of one cycle, exceptions are identified in the outer loop. The corresponding remainders of the selection criteria of these patients with regard to ideal $\mathcal{I}$ are added to the list of generators of ideal $\mathcal{I}$. The polynomials are shown in the third column. The respective patient numbers are shown in the fourth column.}
\end{table}

Here we present the reasons for the choice of relevant true logical statements and for the choice of exceptions.

In the first cycle there have only been two rules (one negative ``$FT(y+1)$'' and one positive ``$Gy(M+x)$'' selection criterion) which can be applied to more than 20 patients and which have less than 5 variables. Thus, we classified them as ``simple rules'' which apply to ``a lot of patients''. This choice is based on our assumption that a small number of parameters should be able to identify sepsis and non-sepsis in the given cohort of patients. The corresponding polynomials of the Gröbner Basis have been $FyTs + FTs$ and $ GMys + Gyxs + GMy + Gyx$ which we added to the list of generators of ideal $\mathcal{J}$. When going through the list of remainders of the selection criteria in Step 4, we only found one polynomial $LyPxTs + LPxTs$ which belongs to a pattern which we only observed once and which is the sum of only two monomials. Adding this polynomial to $\mathcal{I}$ provides a further ``simple'' true logical statement ``$LyPxTs = LPxTs$'' which accounts for all patients, except for patient No. 2237. 
The shared involvement of these molecules in the innate immune response, their induction by similar stimuli, and their connection to TLR/MyD88 signaling pathways make it highly likely that elevated levels of Leptin, TRAIL, Pentraxin-3, and PD-L1 are associated with increased MyD88 expression -- particularly in the context of systemic inflammation such as sepsis. Therefore, we decided to denote patient No. 2237 as an exception.

In the second cycle, the Gröbner Basis of $\mathcal{I}$ changed and the positive and negative rules found addressed more patients compared to the first cycle. We identified two negative rules (selection criterion for non-sepsis patients $T(y+1)(L+M)$ and $E(y+1)(L+M)$) which select more than 30 patients and comprise of ``only'' 4 variables. The corresponding polynomials are shown in Tab.~\ref{tab:results}. In the outer loop we found one polynomial $yxTs+xTs$ which comprises of only 2 monomials and of only 3 variables and ``$s$''. This remainder polynomial was only observed twice (patient No. 127 and patient No. 545). Note that this polynomial is a divider (a generalization) of $LyPxTs+LPxTs$ identified in the first cycle. 

In the last cycle we have only identified one negative rule (the selection criterion $T(y+1)(E+M+P)$) selecting more than 40 patients and ``only'' having 5 variables. The corresponding polynomial $EyTs + MyTs+ yPTs + ETs + MTs + PTs$ of the Gröbner Basis has been added to the list of generators of $\mathcal{J}$. We did not find any ``exception'' anymore, because none of the remainder polynomials had 2 monomials and less than 4 variables plus the variable ``$s$''.

The algorithm terminated and we computed the Gröbner Basis of the ideal $\mathcal{J}$ to figure out the positive and negative rules contained in this basis. They are shown in Tab.~\ref{tab:basis}. The presented selection criteria are very ``restrictive''. 
\begin{table}[hbt]
    \centering
    \begin{tabular}{|l|l|p{6cm}|}
    \hline
    type (pos/neg) & number (sepsis) & selection criterion\cr
    \hline
    neg & 45(1) & T    (y + 1)    (L + M)\cr
    \hline
    neg & 24(0) &   F    T    (y + 1)\cr
    \hline
    neg& 49(0) &   T    (y + 1)    (E + M + P)\cr
    \hline
    pos& 22(22)&   G    y    (M + x) \cr
    \hline
    neg& 33(1) &  E    (y + 1)    (L + M)\cr
    \hline
     pos&  1(1) &   G    y    T    (L  M + L  x + M  x + M)\cr
     neg & 1(1) &   G    T    (y + 1)    (L  M + L  x + M  x + M)\cr
     \hline
    pos & 2(2) &   F    G    y    T    (M + x) \cr
    neg &  3(0) &   F    G    T    (y + 1)    (M + x)\cr
    \hline
    pos & 1(1) &   G    y    T    (E  M + E  x + M  P + M  x + M +
P  x)\cr
    neg &  3(0) &   G    T    (y + 1)    (E  M + E  x + M  P + M  x
+ M + P  x) \cr
    \hline
    pos & 1(1) &  E    G    y    (L  M + L  x + M  x + M)\cr
    neg &  2(1) &   E    G    (y + 1)    (L  M + L  x + M  x + M)\cr
    \hline
    \end{tabular}
    \caption{\label{tab:basis} After the outer loop is finalized, the Gröbner Basis of the ideal $\mathcal{J}$ of ``selected insights'' can be computed. These are all polynomials listed in the second column of Tab.~\ref{tab:results}. Nine polynomials of this basis include a variable ``$s$'' and provide positive and/or negative rules about sepsis. One can see that some negative rules like ``If the patient has a high concentration of ENA-78 and either of Leptin or MMP-8 but not a high concentration of MyD88, then this patient does not have sepsis'' is only valid because of taking out exceptional patients. By this rule $33$ patients are selected, but one of them (No. 2237) has sepsis. The rule is modelled by the selection criterion $E    (y + 1)    (L + M)$ in the fifth row.}
\end{table}

As a follow-up to the results (in the implementation denoted as ``Step 5''), we attempted to further generalize the rules identified in Tab.~\ref{tab:basis}. One can try to further simplify (generalize) the rules of Tab.~\ref{tab:basis} by taking dividers of the selection criteria of negative and positive rules. In this case the number of patients denoted as ``exceptions'' increases. The results of this generalization can be seen in Tab.~\ref{tab:basisfactors}. The negative rule ``$T(y+1)$'' selects 108 patients, but 6 of them have sepsis. The three exceptions from Tab.~\ref{tab:results} are among these six patients. The positive rule ``$Gy$'' selects 67 patents, but only 62 have sepsis. The 5 exceptions do not include the patients No. 2237, 127, or 545. 
\begin{table}[h!]
    \centering
    \begin{tabular}{|l|l|l|p{6cm}|}
    \hline
    type (pos/neg) & number (sepsis) & selection criterion\cr
    \hline
    neg & 108(6) & $T (y + 1)$\cr
    neg & 110(13) & $(y + 1) (E + M + P)$\cr
    neg & 102(12) & $(y + 1) (L + M)$\cr
    neg & 253(39) & $y+1$ \cr
    \hline
    pos & 67(62) & $G y$\cr
    pos & 58(40) & $y(M+x)$\cr
    pos & 50(37) & $G (M + x)$\cr
    pos & 137(98) & $y$\cr
    \hline
    \end{tabular}
    \caption{\label{tab:basisfactors} Positive and negative rules which are generalizations of the rules in Tab.~\ref{tab:basis}.}
\end{table}

\section*{Discussion}

Logical data analysis is based on constructing an ideal $\mathcal{I}$ of selection criteria which select empty sets when applied to the already analyzed patients. We have shown that logical data analysis provides explicit rules which can be used for deriving insights about the relation between protein concentrations and predicting sepsis. It is worth mentioning here, that such rules for a positive classification and for a negative classification differ. A generalized rule in Tab.~\ref{tab:basisfactors} including MyD88 and TRAIL provides a classification as ``non-sepsis'', whereas GLP-1 and MyD88 is used to classify patients as ``sepsis'' cases. Logical data analysis shows how to combine these parameters with further measurements (of Fractalkine, MMP-8, and of Pentraxin-3) in order to get rules with only a few exceptions in Tab.~\ref{tab:basis}.  To include expert knowledge when filtering for relevance is an important step in our method, because some logical connections might be overfittings or artefacts of the concrete dataset.  

Are the found rules reasonable? Tumor necrosis factor-related apoptosis-inducing ligand (TRAIL) and myeloid differentiation primary response 88 (MyD88) are two key molecules involved in immune regulation. Their concentrations in blood plasma have been investigated as potential biomarkers for sepsis diagnosis and prognosis.

\begin{align*}
\text{High TRAIL levels} &\Rightarrow \text{Anti-inflammatory and apoptotic signaling} \\ 
                       & \Rightarrow \text{Reduced systemic inflammation} \\
\text{Low MyD88 levels} &\Rightarrow \text{Attenuated TLR signaling} \\                            &\Rightarrow \text{Lower pro-inflammatory cytokine production}
\end{align*}

Berg et al. \cite{Berg2023} emphasize that TRAIL does not act exclusively in a pro- or anti-inflammatory manner, but rather assumes a balancing role within the complex immunological network of sepsis. TRAIL can induce apoptosis in immune cells, contributing to the regulation of the immune response. Conversely, MyD88 is a central adaptor in Toll-like receptor (TLR) signaling pathways that amplify inflammatory responses. High MyD88 expression has been linked to poor outcomes in sepsis due to excessive cytokine release and organ damage \cite{carlesso2016myd88, mylonas2025myd88}. Therefore, the biomarker profile of ``{\em high TRAIL} and {\em low MyD88} '' ($T(y+1)$) may suggest a non-septic state or a better prognosis, as it indicates a more regulated immune response and reduced risk of cytokine storm.

Glucagon-like peptide-1 (GLP-1) and MyD88 are increasingly recognized as biomarkers associated with sepsis. Elevated levels of both molecules have been linked to poor outcomes in critically ill patients.

\begin{align*}
\text{High GLP-1 levels} &\Rightarrow \text{Metabolic dysregulation and immune modulation} \\
& \Rightarrow \text{Associated with increased mortality} \\
\text{High MyD88 levels} &\Rightarrow \text{Amplified TLR signaling}\\ & \Rightarrow \text{Excessive cytokine release and inflammation}
\end{align*}

GLP-1, traditionally known for its role in glucose homeostasis, has been shown to rise in septic patients, potentially reflecting stress-induced metabolic shifts and immune activation \cite{jacquier2025glp1}. The co-elevation of ``{\em high GLP-1} and {\em high MyD88}'' ($Gy$) may therefore serve as a biomarker signature for sepsis, indicating both metabolic and immunological dysregulation.

The discussed rules seem to have ``only a few exceptions''. They show an important role of MyD88 in discriminating between sepsis and non-sepsis patients in our cohort. However, MyD88 is seemingly not sufficient as a single biomarker according to the number of exceptions (39) when using $y+1$ and $y$ in Tab.~\ref{tab:basisfactors} and thresholds according to Tab.~\ref{tab:extremes}. 

Recent research \cite{BAUER2025106599} has explored the diagnostic potential of proteins including MyD88, GLP-1, and TRAIL in the context of sepsis. The findings indicate that MyD88, Pentraxin-3, PD-L1, and GLP-1 are promising biomarkers, while TRAIL does not offer sufficient diagnostic value. High accuracy in identifying bacterial infections, sepsis, and predicting 30-day mortality was found, amongst others, for MyD88 and GLP-1. The study concludes that combining immunological and metabolic markers -- specifically including MyD88 and GLP-1 -- can significantly enhance the early detection and risk stratification of sepsis which is in accordance with our identified logical rules. In contrast, TRAIL exhibited poor diagnostic performance in \cite{BAUER2025106599} for infection, for sepsis, and for mortality. This is also supported by our logical analysis. We find that high TRAIL combined with low MyD88 concentrations excludes sepsis in many cases in our cohort which could explain why TRAIL is not a single ``positive'' marker for sepsis. A deeper analysis of the polynomials of $\mathcal{I}$ might lead to identify more TRAIL-connected interdependencies in the given set of parameters which are readable from the Gr\"obner Basis. For example, the polynomial $MxT+MxTy$ is element of the Gröbner Basis of $\cal I$. This means that a high level of MMP-8, Pentraxin-3, and TRAIL is always connected to a high concentration of MyD88 in our case study. MyD88 is a central hub \cite{Cheng2020} in the immune response. It acts as a pivotal adaptor protein in innate immunity, influencing key downstream genes such as MMP-8, Pentraxin-3, and potentially TRAIL. 

\subsection*{Limitations}
A prerequisite for applying the presented method is the usefulness of the initial thresholds in Tab.~\ref{tab:extremes}. This approach only works if this choice is meaningful. We plan to extend the method and to include an optimization of these thresholds.

Furthermore, it is important to mention that the choice of relevant true logical statements and exceptions could be different for a different group of experts looking at the polynomials. 

The Gr\"obner Basis is complete. This means, it generates the whole set of true logical statements for the given dataset.  However, the ``simple rules'' may become visible only in the course of the proposed algorithm.

\section*{Acknowledgments}
All the authors have accepted responsibility for the entire content of this manuscript and approved its submission.  MW and KF worked out the logical data analysis method; KK and WB provided the measurements of this study and expert knowledge in the field of sepsis; MR worked on the research data management and on the algebraic workflow description; AG discussed limitations, the role of exceptions, and possible extensions of the method regarding thresholds. The work of MR was supported by the Mathematical Research Data Initiative (MaRDI), part of the National Research Data Infrastructure (NFDI), funded by the German Research Foundation (DFG, project 460135501). KF likes to thank MATH+ (Funded by the Deutsche Forschungsgemeinschaft (DFG, German Research Foundation) under Germany´s Excellence Strategy – The Berlin Mathematics
Research Center MATH+ (EXC-2046/1, EXC-2046/2, project ID: 390685689)).
The authors report no conflicts of interest relevant to this work.

\subsection*{Ethics statement}
The study was approved by the ethics committee of the Charité – Universitätsmedizin Berlin (EA2/044/20 and EA4/167/18), and the study was conducted in accordance with the Declaration of Helsinki (as revised in 2013). All participants were adults. For all patients capable of informed decision-making, written informed consent was obtained prior to participation. Participants received oral and written information about the study, had the opportunity to ask questions, and subsequently signed the consent form. A copy of the signed documents was provided to them. In emergency situations in which patients were temporarily unable to provide consent, the ethics committee allowed the collection and temporary storage of a blood sample until written consent could be obtained from the patient or their legal representative. If written consent was not obtained afterwards, the samples were discarded. Any previously expressed refusal by a patient was respected and led to exclusion.
The study did not include minors, and no waivers of consent beyond the procedure described above were granted. Only those authors which are clinicians at Charité (WB and KK) had access to information that could identify individual participants during or after data collection. On 30/09/2024, the data were previously converted into $\{0,1\}$-pattern. The dataset used for the submitted manuscript does not contain information enabling for subsequent identification of patients. More information is given in Table \ref{tab:ethics}.  

\begin{table}
\centering
\begin{tabular}{|l|l|l|l|l|}
\hline
study & ethics & final vote & recruitment start & recruitment end\\
\hline
HostDx & EA4/167/18 & Jan 15, 2019 &	 Feb 04, 2019 & Sep 04, 2020\\
CASCADE & EA2/044/20 & May 26, 2020 & May 18, 2020 & Aug 31, 2022\\
\hline
\end{tabular}
\caption{\label{tab:ethics} Time-line of the used studies.}
\end{table}

\section*{Supporting material}

Additional material supporting the findings of this study is available from the corresponding author upon reasonable request. The supporting material comprises of the implementation of the algorithm:
The used data set is available on request.
\begin{itemize}
\item[] {\sf Diagnosis\_for\_Sepsis.ipynb} is a Jupyter Notebook based on Julia/OSCAR which implements the presented algorithm. {\sf Diagnosis\_for\_Sepsis.pdf} is a print out of the notebook after the last cycle of the algorithm.
\item[] The used data set of the patient numbers and their $\{0,1\}$-patterns is available on request.
\end{itemize}

\bibliography{sample} 

\end{document}